\documentclass{article}

\usepackage{epsfig}

\usepackage{latexsym}

\begin{document}

\pagestyle{myheadings} \markright{\sc Circular chromatic index
\hfill} \thispagestyle{empty}

\newtheorem{theorem}{Theorem}[section]

\newtheorem{corollary}[theorem]{Corollary}

\newtheorem{guess}[theorem]{Conjecture}

\newtheorem{claim}[theorem]{Claim}

\newtheorem{problem}[theorem]{Problem}

\newtheorem{question}[theorem]{Question}

\newtheorem{lemma}[theorem]{Lemma}

\newtheorem{proposition}[theorem]{Proposition}

\newtheorem{prob}{Problem}

\newcommand{\remark}{\medskip\par\noindent {\bf Remark.~~}}

\newcommand{\pp}{{\it p.}}

\newcommand{\de}{\em}

\newtheorem{example}{Example}

\newenvironment{proof}{{\bf Proof.}}{\hfill\rule{2mm}{2mm}}


\newenvironment{definition}{{\bf Definition.}}{}

\def\mod{ {\ \rm mod \ }}

\def\div{ {\ \rm div \ }}

\newtheorem{prelem}{{\bf Theorem}}

\renewcommand{\theprelem}{{\Alph{prelem}}}

\newenvironment{lem}{\begin{prelem}{\hspace{-0.5

               em}{\bf.}}}{\end{prelem}}

\newtheorem{prelemm}{{\bf Lemma}}

\renewcommand{\theprelemm}{{\Alph{prelemm}}}

\newenvironment{lemm}{\begin{prelemm}{\hspace{-0.5

               em}{\bf.}}}{\end{prelemm}}

\title{Circular chromatic index of graphs of maximum degree $3$}

\author{Peyman Afshani\footnote{Department of Computer Science, University
of Waterloo},
Mahsa Ghandehari\footnote{Department of Industrial Engineering,
Sharif University of Technology},
Mahya Ghandehari\footnote{Department of Mathematics \& Statistics,
Concordia University},
Hamed Hatami\footnote{ Department of Computer Science, University
of Toronto},
\\
Ruzbeh Tusserkani\footnote{Department of Computer Engineering,
Sharif University of Technology}, and
Xuding Zhu\footnote{Department of Applied Mathematics, National
Sun Yat-sen University, Kaohsiung, Taiwan 80424, \newline \ \
e-mail: zhu@math.nsysu.edu.tw} \thanks{This research was partially
supported by the National Science Council under grant
NSC92-2115-M-110-007}}

\maketitle


\begin{abstract}
This paper proves that if $G$ is a graph (parallel edges allowed)
of maximum degree $3$, then $\chi_c'(G) \leq 11/3$ provided that
$G$ does not contain $H_1$ or $H_2$ as a subgraph,  where $H_1$
and $H_2$ are obtained by subdividing one edge of $K_2^3$ (the
graph with three parallel edges between two vertices) and $K_4$,
respectively.  As $\chi_c'(H_1) = \chi_c'(H_2) = 4$, our result
implies that there is no graph $G$ with $11/3 < \chi_c'(G) < 4$.
It also implies that if $G$ is a $2$-edge connected cubic graph,
then $\chi'(G) \le 11/3$.
\end{abstract}

\section{Introduction}

Graphs considered in this paper may have parallel edges but no
loops. Given a graph $G=(V, E)$, and positive integers $p \geq q$,
a $(p, q)$-coloring of $G$ is a mapping $f: V \to \{0, 1, \cdots,
p-1\}$ such that for every edge $e=xy$ of $G$, $q \leq |f(x) -
f(y)| \leq p-q$. The {\em circular chromatic number} $\chi_c(G)$
of $G$ is defined as
$$\chi_c(G) = \inf\{p/q: G
\mbox{\rm \  has a $(p, q)$-coloring} \}.$$
It is known \cite{vince,Zhus} that for any graph $G$, the infimum
in the definition is always attained and
$$\chi(G) - 1 < \chi_c(G) \leq \chi(G).$$
 For a graph $G =(V, E)$, the {\em line graph}
$L(G)$ of $G$ has vertex set $E$, in which $e_1 \sim e_2$, if
$e_1$ and $e_2$ have an end vertex in common. The {\em circular
chromatic index} $\chi_c'(G)$ of $G$ is defined as
$$\chi_c'(G) = \chi_c(L(G)).$$

Recall that the {\em chromatic index} $\chi'(G)$ of $G$ is defined
as $\chi'(G) = \chi(L(G))$. So we have
$$\chi'(G) - 1 < \chi_c'(G) \leq \chi'(G).$$
If $G$ is connected and $\Delta(G) =2$, then $G$ is either a cycle
or a path. This implies that either $\chi_c'(G) = 2$ or
$\chi_c'(G) = 2+\frac{1}{k}$ for some positive integer $k$. Since
graphs $G$ with $\Delta(G) \geq 3$ have $\chi_c'(G) \geq 3$,
`most' of the rational numbers in the interval $(2, 3)$ are not
the circular chromatic index of any graph. The following question
was asked in \cite{Zhus}:

\begin{question}
\label{q1}
 For which rational $r
\geq 3$, there is a graph $G$ with circular chromatic index $r$?
In particular, is it true that for any rational $r \geq 3$, there
is a graph $G$ with $\chi_c'(G) = r$?
\end{question}

If $3<\chi_c'(G)<4$, then $G$ has maximum degree $3$. It is
well-known that the Four Color Theorem is equivalent to the
statement that every  $2$-edge connected cubic planar graph $G$
has $\chi'_c(G) = 3$. For nonplanar $2$-edge connected cubic
graphs, Jaeger \cite{jaeger88} (see also page 197 of \cite{toft})
proposed the following conjecture (Petersen Coloring Conjecture):

\begin{guess}
\label{jaeger} If $G$ is a $2$-edge connected cubic graph, then
one can color the edges of $G$, using the edges of the Petersen
graph as colors, in such a way that any three mutually adjacent
edges of $G$ are colored by three edges that are mutually adjacent
in the Petersen graph.
\end{guess}

Since the Petersen graph has circular chromatic index $11/3$,
Conjecture \ref{jaeger} would imply that every $2$-edge connected
cubic graph $G$ has $\chi'_c(G) \leq 11/3$.  The following two
open problems are proposed in~\cite{Zhus}:

\begin{question}
\label{planar} Prove that if $G$ is a $2$-edge connected cubic
planar graph, then $\chi'_c(G)<4$, without using the Four Color
Theorem.
\end{question}

\begin{question}
\label{bridge} Are there any $2$-edge connected cubic graph $G$
with $\chi'_c(G) =4$?
\end{question}

\begin{figure}
\begin{center}
\includegraphics[width=3in]{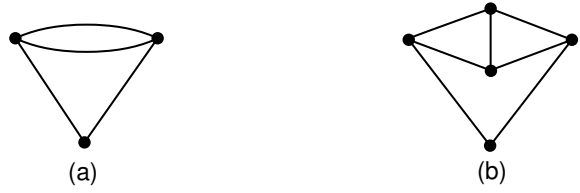}
\end{center}
\caption{\label{Hgraphs}(a): The graph $H_1$, (b): The graph
$H_2$.}
\end{figure}

This paper proves the following result:

\begin{theorem}
\label{main0} Let $H_1$ and $H_2$ be the graphs as shown in
Figure~\ref{Hgraphs}. If $G$ is graph of maximum degree $3$ and
$G$ does not contain $H_1$ or $H_2$ as a subgraph, then
$\chi_c'(G) \leq 11/3$.
\end{theorem}

It is easy to verify that $\chi_c'(H_1) =
\chi_c'(H_2) = 4$. Since graphs $G$ with $\Delta(G) \geq 4$ have
$\chi_c'(G) \geq 4$, we have the following corollary:

\begin{corollary}
\label{gap}
There is no graph $G$ with $11/3 <
\chi_c'(G) < 4$.
\end{corollary}

Corollary \ref{gap} answers the second part of Question~\ref{q1}
in the negative.

To prove Theorem \ref{main0}, it suffices to consider $2$-edge
connected graphs. Indeed, if a graph $G$ is not $2$-edge
connected, say $e$ is a cut edge of $G$, then either $e$ is a
hanging edge, i.e., incident to a degree $1$ vertex, or $e$ is a
cut vertex in $L(G)$. In the latter case, $\chi_c(L(G)) =
\max\{\chi_c(B): B$ is a block of $L(G)\}$. If $e$ is a hanging
edge of $G$, then $e$ has degree at most $2$ in $L(G)$, and hence
any $(11, 3)$-coloring of $L(G)-e$ can be extended to a $(11,
3)$-coloring of $L(G)$. In the remainder of this paper,  we assume
that $G$ is $2$-edge connected and hence has minimum degree at
least $2$. It is easy to see that if $G$ is $2$-edge connected and
has maximum degree at most $3$, then $G$ cannot contain $H_1$ or
$H_2$ as a \emph{proper} subgraph. Therefore Theorem \ref{main0}
is equivalent to the following:

\begin{theorem}
\label{main1}
Suppose $G$ is $2$-edge connected and has maximum degree $3$. If
$G \neq H_1, H_2$, then $\chi'_c(G) \leq 11/3$.
\end{theorem}

Theorem \ref{main1} implies the following corollary, which answers
Questions~\ref{planar} and~\ref{bridge}.

\begin{corollary}
\label{bridgeless}  The circular edge chromatic number of every
$2$-edge connected cubic graph $G$ is less than or equal to
$11/3$.
\end{corollary}

\section{Cubic graphs of girth at least $4$}

The remainder of the paper is devoted to the proof of Theorem
\ref{main1}. In this section, we consider triangle free cubic
graphs. First we prove a lemma which is needed in our proof.

Suppose $c$ is a $k$-coloring of a graph $G=(V, E)$ with colors
$0, 1, \cdots, k-1$.  If $xy$ is an edge of $G$ and $c(y)=c(x)+1
\pmod{k}$, then we say $\vec{xy}$ is a {\em tight arc with respect
to $c$}. Let $A$ be the set of tight arcs, and let $D_c(G)=(V,
A)$, which is a directed graph with vertex set $V$. It is known
\cite{Guichard,Zhus} that if there is a $k$-coloring $c$ of $G$ for which
 $D_c(G)$ is acyclic, then $\chi_c(G)
< k$. The following lemma is a strengthening of this result.

\begin{lemma}
\label{path} Let $c$ be a $k$-coloring of a graph $G$ with colors
$0, 1, \cdots, k-1$,
where $k>2$.
 If $D_c(G)$ is acyclic and each directed path of
$D_c(G)$ contains at most $n$ vertices of color $k-1$, then
$\chi_c(G) \leq k-\frac{1}{n+1}$.
\end{lemma}
\begin{proof}
Let $p=k(n+1)-1$ and $q=n+1$. It suffices to give an $(p,
q)$-coloring for $G$. For each vertex $v$ of $G$, let $l(v)$ be
the maximum number of vertices with color $k-1$ on a directed path
of $D_c(G)$ which ends in $v$, without considering $v$ itself. We
claim that the coloring $c'$ defined as
$$c'(v)=(c(v)q+l(v)) \mod p$$
is a proper $(p, q)$-coloring of $G$. Consider two adjacent
vertices $u$ and $v$. If $2 \leq |c(u)-c(v)| \leq k-2$, then since
both $l(u)$ and $l(v)$ are less than $q$, we have $q \le
|c'(u)-c'(v)|\le p-q$. If $c(u)-c(v)=1$, then $\vec{vu}$ is a
tight arc and hence $l(u) \ge l(v)$. So we have $q \le
|c'(u)-c'(v)|\le p-q$. Finally, if $c(u)=0$ and $c(v)=k-1$, then
$\vec{vu}$ is a tight arc and $l(u) \ge l(v)+1$. Again we have $q
\le |c'(u)-c'(v)|\le p-q$.
\end{proof}

Suppose $c$ is a $k$-edge coloring of $G$ and $e=xy$ is an edge of
$G$. The two arcs $\vec{xy}$ and $\vec{yx}$ are called arcs
corresponding to $e$. We say an arc $\vec{xy}$ is {\em unblocked
with respect to $c$}, if there is a directed walk $W=(e_1, e_2,
\cdots, e_n, e, e'_1, e'_2, \cdots, e'_m)$ in $D_c(L(G))$ such
that (i) $c(e_1) =c(e'_m)=k-1$, and (ii) $e_n=x'x$ and $e'_1=yy'$.
The arc $\vec{xy}$ is {\em blocked with respect to $c$} if no such
directed walk exists. An edge $e=xy$ is said to be {\em blocked in
the direction $x \to y$ with respect to $c$}, if the arc
$\vec{xy}$ is blocked. An edge $e=xy$ is {\em completely blocked
with respect to $c$}, if both arcs $\vec{xy}$ and $\vec{yx}$ are
blocked. Given a partial $k$-edge coloring $c'$ of $G$ (i.e., $c'$
colors a subset of edges of $G$), we say an arc $\vec{xy}$ is
unblocked with respect to $c'$, if $c'$ can be extended to a
$k$-edge coloring $c$ of $G$ such that $\vec{xy}$ is unblocked
with respect to $c$. If no such extension exists, then we say
$\vec{xy}$ is blocked with respect to $c'$. Similarly, we say an
edge $e$ is completely blocked with respect to $c'$, if both arcs
$\vec{xy}$ and $\vec{yx}$ are blocked with respect to $c'$.

\begin{theorem}
\label{main2} If $G$ is a cubic graph of girth at least $4$ and
has a perfect matching, then $\chi'_c(G) \leq 11/3$.
\end{theorem}
\begin{proof}
By Lemma~\ref{path} it suffices to prove that there exists a
$4$-edge coloring $\phi$ of $G$ such that $D_{\phi}(L(G))$ is
acyclic and each directed path of $D_{\phi}(L(G))$ contains at
most two vertices (i.e., two edges of $G$) which are colored by
$3$.

Let $M$ be a perfect matching of $G$. Then $G-M$ is a collection
of cycles. A $4$-edge coloring of $G$ is called a {\em valid
coloring} with respect to $M$, if the following hold:
\begin{itemize}
\item All the $M$-edges (an edge in $M$ is called an
$M$-edge)  are colored by color $0$.\\
\item The edges of any even cycle $C$ of $G-M$ are colored by
colors $1$ and $2$.\\
\item The edges of any odd cycle $C$ of $G-M$ are colored by
colors $1$ and $2$, except one edge which is colored by color $3$.
\end{itemize}

Let $c'$ be a partial $4$-edge coloring of $G$ which can be
extended to a valid $4$-edge coloring of $G$ with respect to $M$.
We are interested in the blocked directions of the $M$-edges with
respect to $c'$. Suppose $e=xy$ is an $M$-edge, and $C$ and $C'$
(not necessarily different)
are cycles of $G-M$ such that $x \in V(C)$ and $y \in V(C')$. If
$\vec{xy}$ is an unblocked arc with respect to $c'$, then we say
$\vec{xy}$ is an {\em input} of $C'$ and an {\em output} of $C$
with respect to $c'$.

Let $C$ be a cycle of $G-M$, and let $c_C$ be the partial edge
coloring of $G$ which is the restriction of a valid coloring $c$
to $M \cup C$. If $C$ is an even cycle, then it is easy to see
that every edge $e \in M$ incident to $C$ is completely blocked
with respect to $c_C$. If $C$ is an odd cycle of $G-M$, then
Figure \ref{blocked} shows the blocked directions of the $M$-edges
incident to $C$ with respect to $c_C$.

\begin{figure}
\begin{center}
\includegraphics[width=3in]{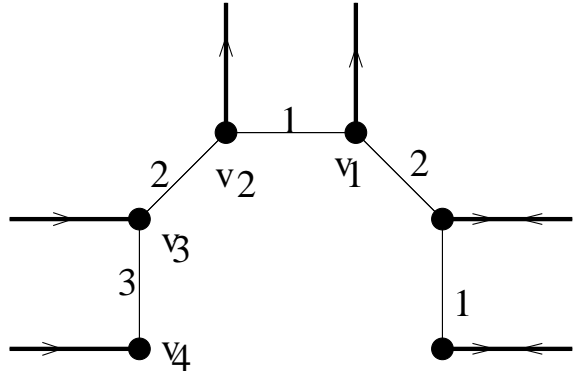}
\end{center}
\caption{\label{blocked}The blocked directions of $M$-edges
incident to $C$ with respect to $c_C$.
}
\end{figure}

In Figure \ref{blocked}, a thick edge indicates an $M$-edge. An
arrow on an $M$-edge indicates a blocked direction of that edge.
An $M$-edge with opposite arrows is completely blocked. Since $G$
has girth at least $4$, the four vertices $v_1, v_2, v_3, v_4$ as
indicated in Figure \ref{blocked} are distinct. Note that an
$M$-edge $e$ incident to  $C$ is completely blocked with respect
to $c_C$, unless $e$ is incident to one of the vertices $v_1, v_2,
v_3, v_4$, which are the vertices on a path whose edges are
colored by colors $1, 2, 3$. So there are at most $4$ $M$-edges
incident to $C$ that are not completely blocked. An $M$-edge
incident to $C$ could be a chord of $C$. If an $M$-edge $e$
incident to $v_1, v_2, v_3, v_4$ is a chord of $C$, then $e$ could
be completely blocked. We will discuss this case later in more
detail. If an $M$-edge $e$ incident to $C$ is not completely
blocked with respect to $c_C$, then exactly one direction of $e$
is blocked.

For a valid $4$-edge coloring $c$ of $G$, let $\phi(c)$ be the
total number of not completely blocked $M$-edges. Let $\psi(c)$ be
the number of not completely blocked $M$-edges that are chords of
cycles of $G-M$.

\begin{claim}
\label{mini}
Suppose $c$ is a valid $4$-edge coloring of $G$ (with respect to
a perfect matching $M$). If $G-M$ has a cycle $C$ which has an
input as well as an output, then there is a valid $4$-edge
coloring $c^*$ of $G$ for which $\phi(c^*)+ \psi(c^*) < \phi(c) +
 \psi(c)$.
\end{claim}
\begin{proof}
Assume $C$ is a cycle of $G-M$ which has an input as well as an
output
with respect to a valid $4$-edge coloring $c$. Then $C$ is an odd
cycle and the $M$-edges incident to $C$ contributes at least $2$
to the summation $\phi(c) + \psi(c)$. We shall construct a valid
$4$-edge coloring $c^*$ of $G$ such that each $M$-edge not
incident to $C$ contributes the same amount to $\phi(c^*)+
\psi(c^*)$ and $\phi(c) +  \psi(c)$. However, the $M$-edges
incident to $C$ contributes at most $1$ to the summation
$\phi(c^*)+ \psi(c^*)$.

Uncolor the edges of $C$ to obtain a partial $4$-edge coloring
$c'$ of $G$. The valid $4$-edge coloring we shall construct is an
extension of $c'$. It is obvious that for any valid $4$-edge
coloring $c^*$ of $G$ which is an extension of $c'$, each
$M$-edge not incident to $C$ contributes the same amount to
$\phi(c^*)+ \psi(c^*)$ and $\phi(c) +  \psi(c)$. So we only need
to make sure that the $M$-edges incident to $C$ contribute at most
$1$ to the summation $\phi(c^*)+ \psi(c^*)$.

First we consider the case that $C$ has no chord. As each
$M$-edge $e$ incident to $C$ is incident to another cycle of
$G-M$, at least one direction of $e$ is blocked with respect to
$c'$. Since $C$ is an odd cycle and $C$ has an input and an
output with respect to $c$, it is easy to see that there are four
consecutive vertices $v_1, v_2, v_3, v_4$ of $C$ such that with
respect to the partial edge coloring $c'$,  the $M$-edges
incident to $v_1, v_2$ have a common blocked direction (i.e.,
either both are blocked in the direction towards $C$ or both are
blocked in the direction away from $C$), and the $M$-edges
incident to $v_3, v_4$ have an opposite blocked direction.
Depending on which directions of the four edges are blocked,
there are four cases as depicted in Figure \ref{block}.

\begin{figure}[ht]
\centerline{\psfig{figure=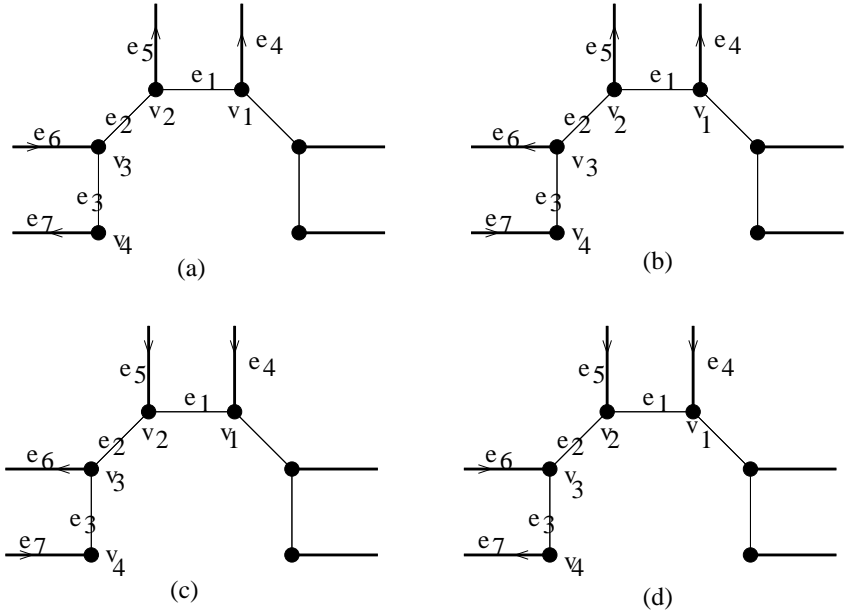}}
\caption{The blocked directions of $M$-edges incident to the uncolored cycle $C$ of $G-M$}
\label{block}
\end{figure}

We use the following convention to interpret Figure \ref{block}
and the figures in the remaining of the paper: An $M$-edge without
an arrow could be completely blocked, or blocked in one direction,
or unblocked in both directions. An $M$-edge with one arrow means
that the indicated direction of that edge is blocked, but the
other direction of that edge could be blocked or unblocked. An
$M$-edge with a pair of opposite arrows means that edge is
completely blocked.

Consider the case indicated in Figure \ref{block} (a) and
\ref{block} (b). We extend $c'$ to a valid $4$-edge coloring $c^*$
of $G$ by letting $c^*(e_1)=3$, $c^*(e_2)=2$, $c^*(e_3)=1$ (the
other edges of $C$ are colored by $1$ and $2$ alternately). It is
easy to verify that in the case indicated in Figure 3(a), $e_7$ is
the only edge which is probably not completely blocked with
respect to $c^*$. In Figure 3(b), $e_6$ is the only edge which is
probably not completely blocked. Thus the $M$-edges incident to
$C$ contributes at most $1$ to the summation $\phi(c^*)+
\psi(c^*)$.

For the cases in Figure 3(c) and 3(d), let $c^*(e_1)=1,
c^*(e_2)=2, c^*(e_3)=3$. Then the $M$-edges incident to $C$
contributes at most $1$ to the summation $\phi(c^*)+\psi(c^*)$.

Next we consider the case that $C$ has a chord.

Since $C$ is an odd cycle, there is an $M$-edge incident to $C$
which is not a chord of $C$. So there is a vertex $v_2$ of $C$
which is incident to a chord of $C$ and a neighbour $v_1$ of $v_2$
in $C$ is not incident to a chord of $C$. Let $v_3, v_4$ be the
vertices of $C$ following $v_1, v_2$ (as shown in Figure
\ref{medge}).

\begin{figure}[ht]
\centerline{\psfig{figure=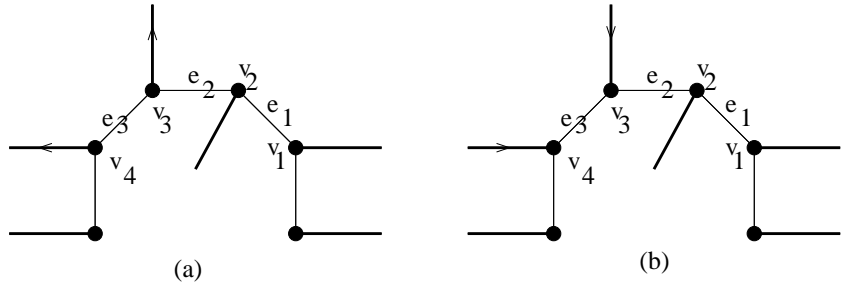}}
\caption{The $M$-edges incident to $v_3, v_4$ have a common blocked direction}
\label{medge}
\end{figure}

Assume the $M$-edges incident to $v_3, v_4$ are not chords of $C$
and have a common blocked direction, as shown in Figure 4(a) or
4(b). In the case as shown in Figure 4(a), extend $c'$ to $c^*$ by
letting $c^*(e_1)=1, c^*(e_2)=2, c^*(e_3)=3$ (and color the other
edges of $C$ alternately by colors $1$ and $2$). In the case as
shown in Figure 4(b), extend $c'$ to $c^*$ by letting
$c^*(e_1)=3, c^*(e_2)=2, c^*(e_3)=1$. In any case, it is easy to
verify that all the chords of $C$ are completely blocked, and
there is at most one $M$-edge incident to $C$ which is not
completely blocked.

Assume the $M$-edges incident to $v_3, v_4$ have opposite blocked
directions or at least one of the $M$-edges incident to $v_3, v_4$
is a chord of $C$. Then depending on which direction of the
$M$-edge incident to $v_1$ is blocked (with respect to $c'$), we
color the edges as in Figure \ref{opposite}.

\begin{figure}[ht]
\centerline{\psfig{figure=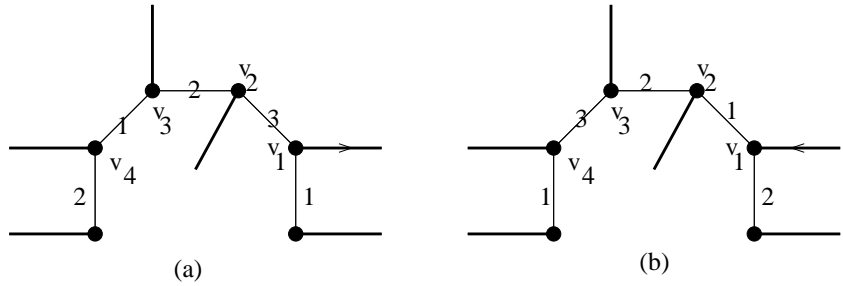}}
\caption{The $M$-edges incident to $v_3, v_4$ have an opposite blocked direction or one of the $M$-edges is a chord.}
\label{opposite}
\end{figure}

In each of the colorings, it is straightforward to verify that
the $M$-edges incident to $C$ contribute at most $1$ to the
summation $\phi(c^*) +
\psi(c^*)$. This completes the proof of Claim \ref{mini}.
\end{proof}

Now we choose a valid $4$-edge coloring $c$ of $G$ such that
$\phi(c)+ \psi(c)$ is minimum. By Claim \ref{mini}, no cycle $C$
of $G-M$ has an input and an output. Since each cycle $C$ of
$G-M$ contains at most one edge of color $3$, it follows that
every directed path of $D_c(L(G))$ contains at most $2$ vertices
(i.e., edges of $G$) with color $3$. By Lemma \ref{path},
$\chi_c(L(G))=\chi'_c(G) \leq 11/3$.
\end{proof}

\begin{corollary}
\label{girth4}
If $G$ is a $2$-edge connected graph of maximum degree $3$ and has
girth at least $4$, then $\chi'_c(G)
\leq 11/3$.
\end{corollary}
\begin{proof}
If $G$ is cubic, then by Petersen Theorem, $G$ has a perfect
matching. Otherwise, take the disjoint union of two copies of $G$,
say $G$ and $G'$. For each degree $2$ vertex $x$ of $G$, connect
$x$ to the corresponding vertex $x'$ in $G'$ by an edge. The
resulting graph $G''$ is cubic (as $G$ has minimum degree $2$) and
is either $2$-edge connected (if $G$ has at least two degree $2$
vertices), or has exactly one cut edge. In any case $G''$ has a
perfect matching (see for example \cite{west}, page 124)
 and has girth at least
$4$. Hence  $\chi'_c(G'') \leq 11/3$ by Theorem \ref{main2}.
\end{proof}

\section{Proof of Theorem \ref{main1}}

We prove Theorem \ref{main1} by induction on the number of edges.
If $|E(G)| =3$, then it is equal to $K_2^3$, and has circular
chromatic index $3$. Assume $|E(G)| \geq 4$ and $G \neq H_1,
H_2$. If $G$ has girth at least $4$, then the conclusion follows
from Theorem \ref{main2}. Thus we assume that $G$ has a pair of
parallel edges or has a triangle.

\begin{figure}
\includegraphics[width=5in]{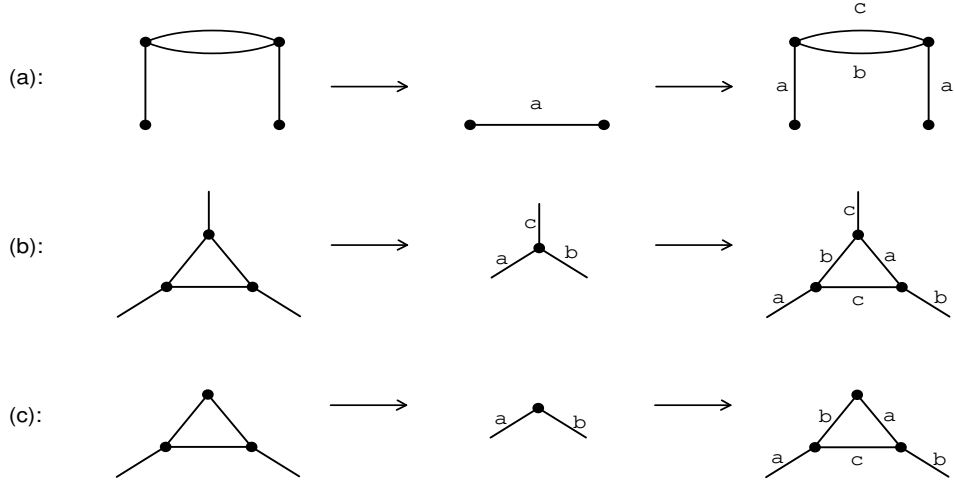}
\caption{\label{operations}(a), (b), and (c) show that how a
$(11/3)$-edge coloring of the new graph leads to a $(11,3)$-edge
coloring of the previous one: (a): In the $(11,3)$-edge coloring
of the main graph $b=(a+3) \mod 11$ and $c=(a+6) \mod 11$, (b):
contracting a triangle with three vertices of degree $3$, (c):
after contracting a triangle with one vertex of degree $2$, we
can always find a color $c$ to complete the $(11,3)$-coloring of
the old graph.}
\end{figure}

\begin{figure}
\begin{center}
\includegraphics[width=5in]{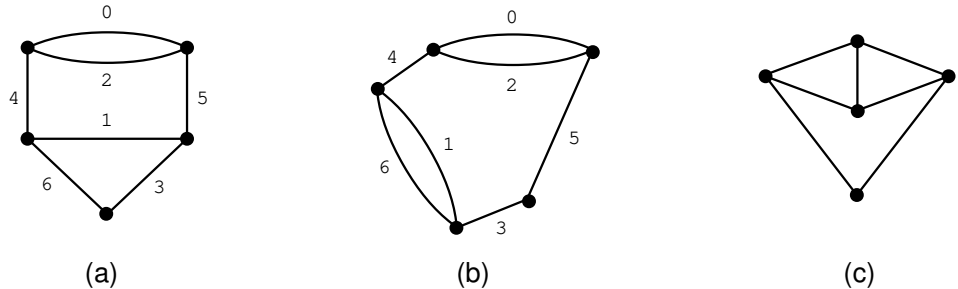}
\end{center}
\caption{\label{h1converts}The graphs that can be converted to
$H_1$ by the ``$\odot$'' operation. For each graph other than
$H_2$ a $(7,2)$-edge coloring is given.}
\end{figure}

\begin{figure}
\begin{center}
\includegraphics[width=5in]{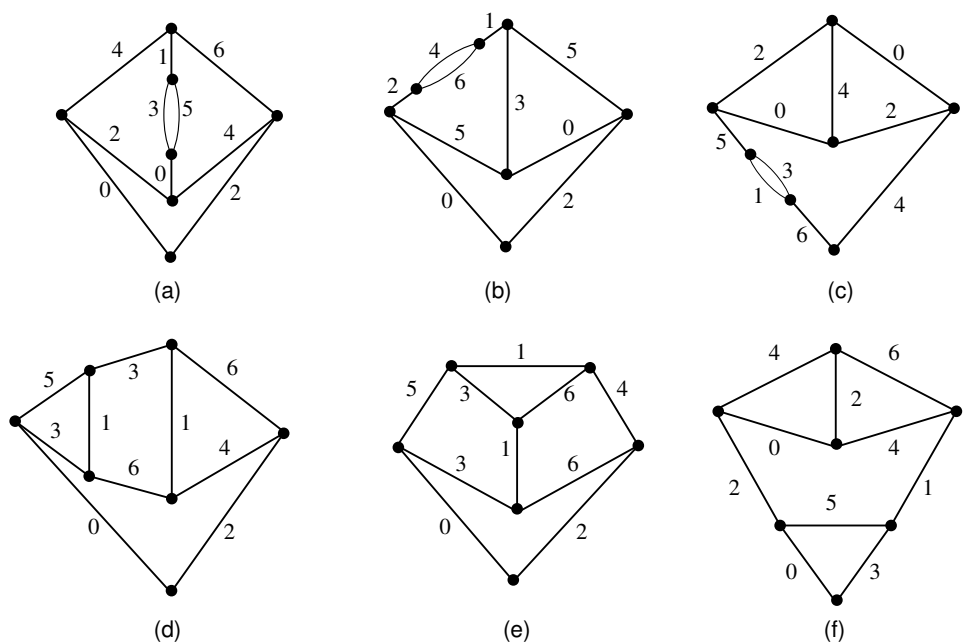}
\end{center}
\caption{\label{h2converts}The graphs that can be converted to
$H_2$ by the ``$\odot$'' operation. For each graph a $(7,2)$-edge
coloring is given.}
\end{figure}

\noindent {\bf Case I:} Suppose there is a pair of parallel edges
between $u$ and $v$. Since $G$ is $2$-edge connected and $G \neq
H_1$, we conclude that $u$ is connected to another vertex $u'$,
$v$ is connected to another vertex $v'$, and $u' \neq v'$. Let
$G\odot uv$ be the graph obtained from $G$ by deleting the two
vertices $u$ and $v$ from $G$ and adding an edge between $u'v'$.
Note that this new edge may cause a multiple edge between $u'$ and
$v'$. If $G \odot uv \not\in \{H_1, H_2\}$, then by induction
hypothesis, $\chi_c'(G \odot uv) \leq 11/3$.
Figure~\ref{operations}(a) illustrates that a $(11, 3)$-coloring
of $L(G \odot uv)$ can be `extended' to a $(11, 3)$-coloring of
$L(G)$. If $G \odot uv \in \{H_1, H_2\}$, then $G$ is one of the
graphs illustrated in Figure~\ref{h1converts} or
Figure~\ref{h2converts}, where a $(7, 2)$-coloring of $L(G)$ is
given.

\noindent {\bf Case II:} Suppose $G$ has a triangle $uvw$. Since
$G$ is $2$-edge connected and $G \neq H_1$, there are no multiple
edges in this triangle. Let $G \odot uvw$ be the graph obtained
from $G$ by contracting the triangle $uvw$ in $G$ to a new vertex.
If $G \odot uvw \not\in \{H_1, H_2\}$, then by induction
hypothesis, $\chi_c'(G \odot uvw) \leq 11/3$.
Figure~\ref{operations}(b,c) illustrates that a $(11, 3)$-coloring
of $L(G \odot uvw)$ can be `extended' to a $(11, 3)$-coloring of
$L(G)$. If $G \odot uvw \in \{H_1, H_2\}$, then $G$ is one of the
graphs illustrated in Figure~\ref{h1converts} or
Figure~\ref{h2converts}, where a $(7, 2)$-coloring of $L(G)$ is
given. So in any case, $\chi_c'(G) \leq 11/3$. This completes the
proof of Theorem~\ref{main1}.

Based on the result in this paper, we propose the following
conjecture:

\begin{guess}
\label{conj}
For any integer $k \geq 2$, there is an $\epsilon > 0$ such that
the open interval $(k-\epsilon, k)$ is a gap for circular
chromatic index of graphs, i.e., no graph $G$ has $k-\epsilon <
\chi'_c(G) < k$.
\end{guess}

If Conjecture \ref{conj} is true, then let $\epsilon_k$ be the
largest real number for which $(k-\epsilon_k, k)$ is a gap for
the circular chromatic index of graphs. The next problem would be
to determine the value of $\epsilon_k$. For $k=2, 3, 4$,
Conjecture
\ref{conj} is true and we know that $\epsilon_2=1, \epsilon_3=1/2$
and $\epsilon_4=1/3$. So a natural guess for $\epsilon_k$ is that
$\epsilon_k = 1/(k-1)$. However, at present time, support for
such a conjecture is still weak. For $k \geq 4$, we do not have
natural candidate graphs $G$ with $\chi'_c(G) = k- 1/(k-1)$.


\bibliographystyle{plain}
\bibliography{circ_edge}

\begin{thebibliography}{1}

\bibitem{Guichard}
D.R. Guichard.
\newblock Acyclic graph coloring and the complexity of the star chromatic
  number.
\newblock {\em J. Graph Theory}, 17:129--134, 1993.

\bibitem{jaeger88}
F.~Jaeger.
\newblock Nowhere-zero flow problems.
\newblock {\em In: L.W.Beineke and Sheehan, editors, Selected Topics in Graph
  Theory}, 3:71--95, 1988.

\bibitem{toft}
T.R. Jensen and B.~Toft.
\newblock {\em Graph Coloring Problems}.
\newblock John Wiley \& Sons, United States of America, 1995.

\bibitem{vince}
A.~Vince.
\newblock Star chromatic number.
\newblock {\em J. Graph Theory}, 12:551--559, 1988.

\bibitem{west}
D.B. West.
\newblock {\em Introduction to Graph Theory}.
\newblock Prentice-Hall, Inc, USA, 2001.
\newblock 2nd Edition.

\bibitem{Zhus}
X.~Zhu.
\newblock Circular chromatic number: a survey.
\newblock {\em Discrete Math.}, 229:371--410, 2001.

\end{thebibliography}

\end{document}